\documentclass[a4paper,11pt]{article}
\usepackage{fullpage}
\makeatletter
\newcommand{\cov}{\mathop{\operator@font cov}}
\makeatother
\begin{document}
\bibliographystyle{plain}
\title{A normal distribution for the disturbance term in regression theory}
\author{Lambros Iossif\\ Athens, Greece}
\date{September 2007}
\maketitle
\begin{abstract}
In regression theory, it is stated that the disturbance term follows the normal distribution 
when the sample size is large. 

In Professor J. Johnston's words: ``In view of the many factors 
involved, an appeal to the Central Limit Theorem would further suggest a normal distribution 
for $u$.''~\cite{johnston63} 

This paper includes an elementary proof that the disturbance term follows the normal distribution 
when $n$ is large.
\end{abstract}

Consider the regression equation
\begin{displaymath}
Y_{i}=\alpha+\beta X_{i}+u_{i}.
\end{displaymath}
The assumptions about the disturbance term $u_i$ are summarized as follows:
\begin{eqnarray*}
E(u_i) &=& 0\\
E(u_i u_j) &=& \left\{ \begin{array}{ll}
                       \sigma^2 & \mbox{when}\;i=j,\, i=1,2,3,\ldots,n\\
                       0        & \mbox{when}\;i\not=j,\, j=1,2,3,\ldots,n
                       \end{array}\right.
\end{eqnarray*}

No assumption is made regarding the distribution of the disturbance term $u_i$.
We wish to show that for large $n$ the term $u_i$ follows the normal distribution
$\mathcal{N}(O,\sigma^2)$ with mean $0$ and variant $\sigma^2$.

We proceed as follows: let
\begin{eqnarray*}
u_i &=&\sum_{k=1}^{n}u_{k} - \sum_{k=1\atop k\not=i}^{m}u_{k}\\
    &=&\frac{n\sigma}{\sqrt{n}}\left[\frac{\displaystyle\sum_{k=1}^{n} \frac{u_k}{n}-O}{%
        \sigma\sqrt{n}}\right] + \frac{\sigma(n-1)}{\sqrt{n-1}}
        \left[ \frac{\displaystyle\sum_{k=1\atop k\not=1}^{n} \frac{-u_k}{n-1}-O}{%
        \sigma/\sqrt{n-1}}
        \right].
\end{eqnarray*}
Let
\begin{displaymath}
x=\frac{\displaystyle\frac{\displaystyle\sum_{k=1}^{n}u_k}{n}-O}{\sigma/\sqrt{n}}
\quad and \quad y=\frac{\displaystyle\sum_{k=1\atop k\not=1}^{n}\frac{-u_k}{n-1}-O
}{\sigma/\sqrt{n-1}}
\end{displaymath}
Then both the random variables $X$ and $Y$ approach the $\mathcal{N}(0, 1)$ distribution
for large $n$. By the Central Limit Theorem:
\begin{displaymath}
f_{X}(x)=\frac{1}{\sqrt{2\pi}}\exp\left(-\frac{x^2}{2}\right)\quad\mbox{and}\quad
f_{Y}(y)=\frac{1}{\sqrt{2\pi}}\exp\left(-\frac{y^2}{2}\right)
\end{displaymath}
for large $n$. The joint density $f_{XY}(x, y)$ is completely determined if we
can calculate the correlation $\rho_{XY}$ between $X$ and $Y$. Thus we have: since
$\sigma_{X} = \sigma_{Y} = 1$
\begin{eqnarray*}
\rho_{XY} &=& \frac{\cov(X,Y)}{\sigma_{X}\sigma_{Y}}\\
          &=& E\left[\left(\frac{\displaystyle\frac{1}{n}\sum_{k=1}^{n}u_{k}}{\sigma/\sqrt{n}} -
                 E\left(\frac{\displaystyle\sum_{k=1}^{n}\frac{u_{k}}{n}}{\sigma/\sqrt{n}}\right)\right)
                  \left(\frac{\displaystyle\sum_{k=1\atop k\not=i}^{n}\frac{-u_k}{n-1}}{\sigma/\sqrt{n-1}}-E
                  \left(\frac{\displaystyle\sum_{k=1\atop k\not=i}^{n}\frac{-u_k}{n-1}}{\sigma/\sqrt{n-1}}
                  \right)\right)\right]\\
          &=& E\left[\left(\frac{\displaystyle\sum_{k=1}^{n}\frac{u_k}{n}}{\sigma/n}\right)
                     \left(\frac{\displaystyle\sum_{k=1\atop k\not=i}^{n}\frac{-u_k}{n-1}}{\sigma/\sqrt{n-1}}
                     \right)\right]\\
          &=& E\left[\left(\frac{\displaystyle\sum_{k=1}^{n}\frac{u_k}{n}}{\sigma/\sqrt{n}}\right)
                \left(\frac{\displaystyle\sum_{k=1}^{n}\left(\frac{-u_k}{n-1}+\frac{u_k}{n-1}\right)}{%
                \sigma/\sqrt{n-1}}\right)\right]\\
          &=& \frac{1}{\sigma^2/\sqrt{n(n-1)}} E\left[\frac{\displaystyle-\left(\sum_{k=1}^{n} u_{k} \right)^2}{%
              n(n-1)}+\frac{\displaystyle u_{i}\sum_{k=1}^{n} u_{k}}{n(n-1)}\right]\\
          &=& \frac{\sqrt{n(n-1)}}{\sigma^2}E\left[\frac{\displaystyle-\sum_{i=1}^{n}\sum_{j=1}^{n}u_i u_j}{%
              n(n-1)} + \frac{\displaystyle u_i \sum_{k=1}^{n}u_k}{n(n-1)}\right]\\
          &=& \frac{\sqrt{n(n-1)}}{\sigma^2}E\left[\frac{\displaystyle-\sum_{i=1}^{n}\sum_{j=1\atop j\not=i}^{n}
              E(u_i u_j)}{n(n-1)} - \frac{\displaystyle\sum_{i=1}^{n}E(u^2_i)}{n(n-1)} +
              \frac{\displaystyle\sum_{k=1\atop k\not=i}^{n}E(u_i u_k)}{n(n-1)} +
              \frac{E(u^2_i)}{n(n-1)}\right]\\
         &=& \frac{\sqrt{n(n-1)}}{\sigma^2}\left[\frac{-0-n\sigma^2 + 0 - \sigma^2}{n(n-1)}\right]\\
         &=&  \frac{-(n-1)\sqrt{n(n-1)}}{n(n-1)}\\
         &=&  \frac{-\sqrt{n(n-1)}}{n}\\
         &=&  \rho_{XY}. 
\end{eqnarray*}
Thus the joint density $f_{XY}$ of $X$ and $Y$ is the bivariate normal density
function:
\begin{displaymath}
f_{XY}(x,y)=\frac{1}{2\pi\sqrt{1-\rho^2}}\exp\left\{-\frac{1}{2(1-\rho^2)}\left[x^2-2\rho x y + y^2 
\right]\right\}
\end{displaymath}
with $\sigma_{x}=1$, $\sigma_{y}=1$, $\mu_{x}=0$, $\mu_{y}=0$, and $\rho=\frac{-\sqrt{n(n-1)}}{n}$ is
\begin{displaymath}
f_{XY}(x,y)=\frac{\sqrt{n}}{2\pi}\exp\left\{-\frac{1}{2}\left( nx^2 + 2\sqrt{n(n-1)}xy+ny^2\right)\right\}
\end{displaymath}
for large $n$.

    In order to find the density of $u_i$ for large $n$ we consider the transformation 
(see~\cite[p.~204]{mood74}):
\begin{eqnarray*}
u_i &=& \frac{n\sigma}{\sqrt{n}}x+\frac{\sigma(n-1)}{\sqrt{n-1}}y=\sqrt{n}\sigma x+\sigma y\sqrt{n-1}\\
v   &=& y
\end{eqnarray*}
Let the above transformation be a one-to-one transformation of the $xy$ plane
onto the $uv$ plane with inverse transformation given by:
\begin{eqnarray*}
x &=& \frac{u_i}{\sigma\sqrt{n}}-\frac{\sqrt{n-1}}{\sqrt{n}}v\\
y &=& v.
\end{eqnarray*}
The Jacobian of the above transformation is:
\begin{displaymath}
J=\left| \begin{array}{cc}
\displaystyle\frac{\partial x}{\partial u_i} & \displaystyle\frac{\partial x}{\partial v}\\
\displaystyle\frac{\partial y}{\partial u_i} & \displaystyle\frac{\partial y}{\partial v}
\end{array}\right|=
\left| \begin{array}{cc}
      \displaystyle\frac{1}{\sigma\sqrt{n}} & \displaystyle-\frac{\sqrt{n-1}}{\sqrt{n}}\\
       0                       & 1
      \end{array}\right|=\frac{1}{\sigma\sqrt{n}}.
\end{displaymath}
We note in passing that the partial derivatives $\partial x/\partial u_i$, 
$\partial x/\partial u$, $\partial y/\partial u_i$, and
$\partial y/\partial u$ are all continuous functions of $u_i$ and $u$ as they are constant.

We are interested in the absolute value of the Jacobian above:
\begin{displaymath}
|J|=\frac{1}{\sigma\sqrt{n}}.
\end{displaymath}
Then the joint density of $u_i$ and $v$ is
\begin{eqnarray*}
f(u_i, v) &=& \frac{\sqrt{n}}{2\pi\sigma\sqrt{n}}\exp\left[-\frac{1}{2}\left\{n\left(
               \frac{u_i}{\sigma\sqrt{n}}-\frac{\sqrt{n-1}}{\sqrt{n}}v\right)^{2}+
               2\sqrt{n(n-1)}\left(\frac{u_i}{\sigma\sqrt{n}}-\frac{\sqrt{n-1}}{\sqrt{n}}v\right)v+nv^{2}
               \right\}\right]\\
         &=& \frac{1}{2\pi\sigma}\exp\left[-\frac{1}{2}\left\{\frac{u_i^2}{\sigma\sqrt{n}}+(n-1)
              v^2 +\frac{2u_i v\sqrt{n-1}}{\sigma}-\frac{2u_i v\sqrt{n-1}}{\sigma}-2(n-1)v^2+nv^2\right\}
              \right]\\
         &=& \frac{1}{\sqrt{2\pi}\sigma\sqrt{2\pi}}\exp\left\{-\frac{1}{2}\left[ 
             \frac{u^2_i}{\sigma\sqrt{n}}+nv^2 + v^2 - 2nv^2 + nv^2\right] \right\}\\
         &=& \frac{1}{\sqrt{2\pi}\sigma}\frac{1}{\sqrt{2\pi}}\exp\left\{
             -\frac{u_i^2}{2\sigma^2}\right\}\exp\left\{-\frac{v^2}{2} \right\}.
\end{eqnarray*}
Thus
\begin{displaymath}
f(u_i)=\frac{1}{\sqrt{2\pi}\sigma}\exp\left\{-\frac{u_i^2}{2\sigma^2} \right\}
\int^{+\infty}_{-\infty}\frac{1}{\sqrt{2\pi}}\exp\left\{\frac{-v^2}{2} \right\}dv,
\end{displaymath}
when $v$ is large.

From our study of the bivariate normal random variable $u$ with density
\begin{displaymath}
f(v)=\frac{1}{\sqrt{2\pi}}\exp\left\{-\frac{u^2}{2}\right\}
\end{displaymath}
the integral above is unity i.e.:
\begin{displaymath}
\int^{+\infty}_{-\infty}\frac{1}{\sqrt{2\pi}}\exp\left\{\frac{-v^2}{2} \right\}dv=1.
\end{displaymath}
So that the density of $u_i$ alone is given by
\begin{displaymath}
f(u_i)=\frac{1}{\sqrt{2\pi}\sigma}\exp\left\{-\frac{1}{2}\frac{u_i^2}{\sigma^2} \right\}=
\mathcal{N}(0,\sigma^2)
\end{displaymath}
with $u$ a normal distribution mean $0$ variance $\sigma^2$.

   This is what we set out to prove when $n$ is large.

   Q.E.D.

\paragraph{Note} We observe that $u_i$ and $v$ are independent as their joint density factors
out as a product of $u_i$ and $v$ alone. This is expected since $u_i$ and
\begin{displaymath}
v=y=\sum_{k=1\atop k\not=i}^{n}\frac{-u_k}{n-1}\Bigm/\frac{\sigma}{\sqrt{n-1}}
\end{displaymath}
are clearly independent.
\nocite{*}
\bibliography{iossif}
\end{document}